\documentclass{amsart}
\usepackage{amssymb,stmaryrd}
\usepackage{amsfonts}
\usepackage{amstext}
\usepackage{geometry}                % See geometry.pdf to learn the layout options. There are lots.
\usepackage{graphicx}    %% packages
\usepackage{scalefnt}
\usepackage{MnSymbol}
\usepackage[all]{xypic}
\usepackage{slashed}
\usepackage{extarrows}
\usepackage{stmaryrd}
\usepackage{lettrine}
\usepackage{booktabs}
\usepackage{paralist}
\usepackage{subfig}
\usepackage{color}

\newtheorem{lemma}{Lemma}[section] 
\renewcommand{\imath}{\mathrm{i}}

\newcommand{\C}{\mathbb{C}}
\newcommand{\R}{\mathbb{R}}

\newcommand{\Z}{\mathbb{Z}}

\newcommand{\del}{\partial}

\newcommand{\ev}{\mathrm{ev}}

\newcommand{\extd}{\mathrm{d}}

\newcommand{\tens}{\mathop{{\otimes}}}

\newcommand{\id}{\mathrm{id}}

\newcommand{\cX}{\mathfrak{X}}

\allowdisplaybreaks
\begin{document}
\title{Quantum geodesic flow on the integer lattice line}

\author{Edwin Beggs and Shahn Majid}


\address{
Department of Mathematics, Bay Campus, Swansea University, SA1 8EN, UK\\
Queen Mary University of London, School of Mathematical Sciences, London E1 4NS, UK\\}

\email{e.j.beggs@swansea.ac.uk,   s.majid@qmul.ac.uk}

\maketitle

\begin{abstract}
 We use a recent formalism of quantum geodesics in noncommutative geometry to construct geodesic flow on the infinite chain $\cdots\bullet$--$\bullet$--$\bullet\cdots$.  We find that noncommutative effects due to the discretisation of the line cause an initially real geodesic flow amplitude $\psi$ (for which the density is $|\psi|^2$) to become complex. This has been noted also for other quantum geometries and suggests that the complex nature of the wave function in quantum mechanics (and the interference effects that follow) may have its origin in a quantum/discrete nature of spacetime at the Planck scale.
\end{abstract}

\section{Introduction}

In recent works\cite{Beg:geo,BegMa:geo,BegMa:cur}, we have introduced a radically new way of thinking about  geodesics even on a classical manifold $M$ (which then extends to noncommutative or `quantum' Riemannian geometry). The idea is to think of not one geodesic at a time but a flow of geodesics much like in fluid dynamics where each particle moves along a geodesic. The tangent vectors to these geodesics form a geodesic-time dependent vector field $X(s)$, where $s$ is geodesic time,  and it turns out that these obey a simple {\em geodesic velocity equation} 
\begin{equation}\label{claveleq} \dot X+\nabla_XX=0\end{equation}
where dot is with respect to $s$.  Next, instead of an actual (evolving) fluid particle density $\rho(s)$, we have an evolving wave function $\psi(s)$ or amplitude with $|\psi(s)|^2=\rho(s)$, and we solve for $\psi$ by the {\em amplitude flow} equation
\begin{equation}\label{claamp} \dot \psi+  \psi \kappa + X(\extd \psi)=0,\quad \kappa=\frac{1}{2}{\rm div}(X).\end{equation}
If one considers bump functions, then classically this reproduces a bump travelling with velocity $X(s)$ evaluated at the bump, i.e. a classical geodesic as expected. If $\psi$ is real-valued then this is not really different from working with $\rho$, i.e. a fluid flow language. But when $\psi$ is complex-valued then there are possible interference effects as normally associated with quantum mechanics. Indeed, we make a quantum mechanics-like interpretation with respect to geodesic time $s$ and (\ref{claamp}) in the role of Schroedinger's equation.

\begin{itemize}
\item This way of thinking of geodesics rips apart the usual concept of geodesics as points with tangent vectors and {\em puts them back in reverse order}, now starting with a velocity field and only afterwards introducing particle densities.

\item Geodesic time $s$ is no longer arclength or proper time of just one particle but the time experienced by an observer who watches all the particles evolving together. It provides a new coordinate for the physical system.

\item If $M$ is spacetime then $\psi(x,t)$  is  (an evolving) spacetime amplitude with $|\psi(x,t)|^2$ the probability density to find a particle at $(x,t)$ at time $s$, which is unfamiliar. However, if the system is $t$-independent then one can focus on wave functions of a fixed frequency in $t$, then the spatial evolution looks more like quantum mechanics with respect to $s$ cf. \cite{BegMa:geo}.

\item The Ricci curvature $R_{\mu\nu}$ enters naturally via the {\em convected derivative} ${D\over D s} f={\del\over\del s} + X(\extd f)$ of a function (the rate of change seen moving with the fluid flow). Then \cite{BegMa:cur}
\[ {D\over Ds}{\rm div}(X)=-(\nabla_\mu X^\nu)(\nabla_\nu X^\mu) - R_{\mu\nu}X^\mu X^\nu\]
\end{itemize}

This new way of thinking about particles in GR is more field-theoretic (starting with the velocity field $X$ as fundamental) but takes  getting used. Here, quantum geodesics on classical Minkowski spacetime and on the bicrossproduct model quantum spacetime $[x,t]=\imath\lambda_P x_i$ in \cite{MaRue} (where $\lambda_P$ is the Planck scale) are studied in \cite{LiuMa}. One of the results there is that a `point' in spacetime modelled as bump with width $\sigma$ moving with respect to $s$ gets a correction to its effective velocity of order $(\lambda_P/\sigma)^2$, i.e. there can be no such thing as a truly point particle {\em due to quantum gravity effects} if these are modelled by quantum spacetime.

Another phenomenon in \cite{LiuMa}, which we support now with new and direct calculations on the integer lattice, is that an initially real $\psi$ {\em does not stay real} when the geometry is quantum. Quantum geometry here means that we replace $C^\infty(M)$ by a possibly noncommutative $*$-algebra $A$. We replace differential forms on $M$ by an $A$-bimodule $\Omega^1$ etc., see below.  We also replace the measure on $M$ with respect to which the probabilistic interpretation of $\psi$ is defined by a reference  {\em state} or positive linear functional $\int:A\to \C$. In this case,  $\psi(s)\in A$ is an element of a $*$-algebra and $|\psi(s)|^2$ is replaced by the positive operator $\psi(s)^*\psi(s)\ge 0$ but nevertheless behaves like a probability density when used to compute expectation values via $\int$. A reality condition on the geodesic velocity field $X$   ensures that the flow is unitary in the sense of preserving $\int \psi^*\psi$.

As well as being applicable to actual quantum mechanics (where the usual Schroedinger equation can be seen as a certain quantum geodesic flow)\cite{BegMa:geo} and to relativistic quantum mechanics, the theory also applies to discrete geometry such as the lattice line in the present work. The algebra of functions $A$ is still commutative, namely functions on the vertices, but the space of differentials $\Omega^1$ is intrinsically noncommutative, being spanned as a vector space  by the arrows $\{x\to y\}$ of the graph. Here, the product on either side by a function $f$ on the vertex set is
\[ f.(x\to y)=f(x)x\to y\quad  \ne\quad f(y)x\to y=(x\to y).f;\]
once we allow such noncommutativity between differentials and functions then discrete geometry sits very naturally as a special case of noncommutative differential geometry. 

Section~\ref{secpre} of the paper outlines the tools of {\em quantum Riemannian geometry} (QRG) as developed in works of the authors and others, see the text \cite{BegMa}, used to write down the quantum versions of (\ref{claveleq}) and (\ref{claamp}). Section~\ref{secZ} then shows how the theory works on the integer line $\Z$, for which the unique $*$-preserving QLC was found in \cite{ArgMa1}. We show by solving the flow equations that, even for the flat lattice where all edge lengths are the same, if $\psi$ starts off real then it evolves into a complex wave function due to the noncommutative nature of the differential geometry of $\Z$. This is also known for geodesics on $A=M_2(\C)$ in \cite{BegMa:cur}, but here it  arises out of the discreteness. This  suggests that the set-up of quantum mechanics with complex wave functions could itself have its origin in quantum gravity leading to spacetime being better modelled\cite{Sny,Ma:pla} as discreteness/noncommutative at the Planck scale.

We note that the methods of QRG used here are very different from the approach to noncommutative geometry of Connes\cite{Con} and others coming out of operator theory and cyclic cohomology. The more constructive QRG approach grew out of \cite{BegMa:rie, BegMa:gra} and more broadly including experience with quantum groups. It was applied in Physics to models of quantum gravity and Beckenstein-Hawking radiation in \cite{Ma:sq, Ma:haw, ArgMa2,LirMa} and to Kaluza-Klein models for elementary particles in \cite{ArgMa4,LiuMa2}. Here \cite{ArgMa4} also provides an introduction for Physicists to the QRG formalism. Quantum geodesics, however, link back to operator theory via the KSGNS construction and KK-theory as a motivation in \cite{Beg:geo}, suggesting a new direction for convergence of these approaches. 
 
 \section{Recap of the formalism}\label{secpre}

We recap elements of QRG from \cite{BegMa} and then of quantum geodesics from \cite{Beg:geo,BegMa:geo, BegMa:cur}. 

\subsection{Quantum Riemannian Geometry} We start with a $*$-algebra $A$, which will play the role of our `coordinate algebra' but could be noncommutative. We require this to have   a differential structure in the form of an $A$-bimodule $\Omega^1$ of differential forms equipped with a map $\extd:A\to \Omega^1$ obeying the Leibniz rule 
\[ \extd(ab)=a\extd b+ (\extd a)b\]
and such that $\Omega^1$ is spanned by elements of the form $a\extd b$ for $a,b\in A$. This can be extended to a full differential graded `exterior algebra' $\Omega$ although not uniquely (there is a unique maximal extension). We assume that at least $\Omega^2$ has been specified  and that  $*$ extends at least to $\Omega^1$ as a graded-involution (i.e., with a minus sign on swapping odd degrees) and commutes with $\extd$. 

Next, a metric means for us an element $\mathfrak{g}\in\Omega^1\tens_A\Omega^1$ which is invertible in the sense of a bimodule map $(\ ,\ ):\Omega^1\tens_A\Omega^1\to A$ obeyng the usual requirements as inverse to $\mathfrak g$, which  forces $\mathfrak g$  to be central\cite{BegMa:gra}. There are ways to go beyond, for example starting with $(\ ,\ )$ or with a metric as a map from $\Omega^1$ to vector fields, as in the hermitian version of the theory\cite[Chap.~8.5]{BegMa}.

Next, a QLC or {\em quantum Levi-Civita connection} is a left bimodule connection $(\nabla,\sigma)$ on $\Omega^1$ which is metric compatible and torsion free in the sense in the sense 
\begin{equation}\label{nablag} \nabla \mathfrak{g}:=(\nabla\tens\id+(\sigma\tens\id)(\id\tens\nabla))\mathfrak g=0,\quad   T_\nabla:= \wedge\nabla-\extd=0\end{equation}
for a left bimodule connection $\nabla:\Omega^1\to \Omega^1\tens_A\Omega^1$ characterised by\cite{DVM,Mou}
\[ \nabla(a.\omega)=\extd a\tens \omega+ a.\nabla\omega,\quad \nabla(\omega.a)=\sigma(\omega\tens\extd a)+(\nabla\omega).a,\]
where the `generalised braiding' $\sigma:\Omega^1\tens_A\Omega^1\to \Omega^1\tens_A\Omega^1$ is assumed to exist and is uniquely determined by the second equation. The thinking behind these formulae is that the left factor in the output of $\nabla$ would classically evaluate against a vector field to provide a covariant derivative along it. There is an analogous theory with right bimodule connections where the right hand factor of the output would be evaluated against a vector field to recover classical geometry.  In terms of $(\ ,\ )$, metric compatibility amounts to 
\begin{equation}\label{metnabla} \extd(\omega,\eta)=(\id\tens(\ ,\ ))(\nabla\omega\tens\eta+(\sigma\tens\id)(\omega\tens\nabla\eta)),\end{equation}
see \cite[Lemma~8.4]{BegMa} for the equivalence with (\ref{nablag}). 

We also require `reality' of the metric and for  $\nabla$ to be $*$-preserving\cite{BegMa:rie,BegMa} in the sense
\[ \mathfrak{g}^\dagger=\mathfrak{g},\quad \sigma\,\dag\,\nabla(\xi^*)=\nabla\xi,\quad  \dagger:={\rm flip}(*\tens *).\]

Finally, we will need $\cX={}_A\hom(\Omega^1,A)$, the space of left quantum vector fields defined as left $A$-module maps  $X: \Omega^1\to A$. This has a bimodule structure 
\[ (a.X.b)(\omega)= (X(\omega.a))b\]
for all $\omega\in \Omega^1$ and $a,b\in A$, and inherits from $\nabla$ a dual right connection  
\[\nabla_\cX:\cX\to \cX\tens_A\Omega^1,\quad \sigma_\cX: \Omega^1\tens_A\cX\to \cX\tens_A\Omega^1.\]
characterised by 
\[  \extd(\ev(\omega\tens X))=(\id\tens\ev)(\nabla \omega \tens X)+(\ev\tens\id)(\omega\tens\nabla_\cX X) \]
for all $\omega\in \Omega^1$ and $X\in \cX$.  We use this to define the divergence of a vector field by
\[ {\rm div}:=\ev\circ\sigma_\cX^{-1}\nabla_\cX,\quad \ev:\Omega^1\tens_A\cX\to A,\quad \ev(\omega\tens_A X)=X(\omega)\]
If $\nabla$ is metric compatible then so is $\nabla_\cX$ in a suitable sense and then one can prove that
\[ {\rm div}(\mathfrak{g}_2(\omega))=(\ ,\ )\nabla \omega;\quad \mathfrak{g}_2:\Omega^1\to \cX,\quad \mathfrak{g}_2(\omega)=(\ ,\omega)\]
for all $\omega\in \Omega^1$, as would be usual in GR  (i.e., one can compute divergences as the codifferential of the 1-form corresponding via the metric to a vector field). The bottom line is, all the basic tools of Riemannian Geometry extend to this set-up provided one is rather careful. 

%Finally, a calculus is called {\em inner} if $\extd=[\theta,\ ]$ for $\theta\in \Omega^1$ with $\theta^*=-\theta$. Then\cite{Ma:gra}\cite[Thm~8.11]{BegMa} a left bimodule connection has the specific form
%\begin{equation}\label{nablainner} \nabla \omega=\theta\tens\omega- \sigma(\omega\tens\theta)+\alpha(\omega)\end{equation}
%for some bimodule maps $\sigma:\Omega^1\tens_A\Omega^1\to \Omega^1\tens_A\Omega^1$ and $\alpha:\Omega^1\to \Omega^1\tens_A\Omega^1$. 

\subsection{Quantum geodesics}\label{secgeo}

We are now ready for quantum geodesics. Although the theory is more general, we will focus on `wave functions' in $\psi\in E=C^\infty(\R, A)$ where $s\in\R$ will be the `geodesic time' parameter. We assume that $A$ is equipped with enough structure for such functions (or some variant of them) to make sense. Likewise, we let $X_s$ be a time-dependent left vector field on $A$ and $\kappa_s$ another time-dependent element of $A$. These obey the {\em geodesic velocity equations} if
\begin{equation}\label{veleqX} \dot X_s +[X_s,\kappa_s]+(\id\tens X_s)\nabla_\cX(X_s)=0\end{equation}
where dot means differential with respect to $s$. 

Next, we require $\int: A\to \C$ to be a non-degenerate positive linear functional (so we can think of it as a probability measure if normalisable so that $\int 1=1$) and define  ${\rm div}_{\int}(X)$ of a vector field by
\begin{equation}\label{divint} \int a\,  {\rm div}_{\int}(X)+ \int X(\extd a)=0\end{equation} 
for all $a\in A$. We  require that the geodesic velocity field and $\kappa$ at each $s$ to obey the {\em unitarity conditions}  
 \begin{equation}\label{unitarity} \int \kappa^* a+a\kappa+X(\extd a)=0,\quad \int X(\omega^*)=\int X(\omega)^*\end{equation}
for all $a\in A$ and all $\omega\in\Omega^1$. Note that if the second of (\ref{unitarity}) applies then one says that $X$ is {\em real with respect to $\int$} and  then we can canonically solve the first equation in the pair by 
\[ \kappa={1\over 2}{\rm div}_{\int}(X),\]
 see after \cite[Def.~4.5]{BegMa:cur}. It is not automatic that if $X_s$ is initially real with respect to $\int$ that it necessarily remains so under the geodesic velocity equation, we have to impose this as a further {\em improved auxiliary equation} obtained as the difference between (\ref{veleqX}) and its conjugate under the reality assumption. 

Given this geodesic velocity vector field, we then require the {\em amplitude flow equation} 
\[ \dot\psi= -\psi\kappa_s-X_s(\extd \psi),\]
where $\extd$ acts on $\psi(s)\in A$ and dot is with respect to $s$ as before. The above conditions ensure that $\int\psi^*\psi$ is constant in time, which is needed for a probabilistic interpretation. 

Finally,  the above works for any $\int$ but in Riemannian geometry we would normally use a particular measure defined by $\sqrt{|g|}$ characterised as vanishing on a total divergence. We do the same and say that  $\int$ is {\em divergence compatible} (with $\nabla$) if  
\[ \int{\rm div} X=0\]
for all $X\in \cX$. This is equivalent to saying that ${\rm div}_{\int}={\rm div}$. If it is also (as it will be in our case) a trace so that $\int (ab)=\int (ba)$ then as a special case of the theory in \cite{BegMa:cur}, $\cX$ acquires a $*$ operation characterised by
\[
 X^*(\omega)=\big(\ev\circ\sigma_\cX^{-1}(X\tens\omega^*)\big)^*
\]
for all $\omega\in\Omega^1$ and such that our previous $X$ being real with respect to $\int$ now appears as $X^*=X$.

\section{Quantum geodesics on $\Z$}\label{secZ}

We illustrate the above machinery on the integer lattice. Here $A=C(\Z)$ denotes all (or some class of) functions on the vertices of the lattice line labelled by integers $\Z$. These are viewed as a graph with arrows in both directions $i\to i+1$ and $i\leftarrow i+1$ for all $i$. In fact this is a Cayley graph and we can formally define $e^+$ as the sum of all increasing arrows and $e^-$ as the sum of all decreasing arrows. We do not need to worry about this sum and can just take $e^\pm$ as, by definition, a basis of $\Omega^1$ over $A$, i.e. every 1-form is of the form $\omega=\omega_+ e^++ \omega_- e^i$ for $\omega_\pm\in A$. We just need not know that the bimodule commutation relations and exterior derivative are
\[ e^\pm f= R_\pm(f) e^\pm,\quad \extd f= (\del_\pm f)e^\pm;\quad R_\pm(f)(i)=f(i\pm 1),\quad \del_\pm=R_\pm-\id.\]
So the $R_\pm$ are just lattice shift up and down and $\del_\pm$ are just the usual lattice finite difference. 

A metric $\mathfrak{g}\in\Omega^1\tens_A\Omega^1$ is then just an element of the form
\[ \mathfrak{g}= g_+ e^+\tens e^- + g_- e^-\tens e^+;\quad g_\pm \in \R\setminus\{0\}\]
where the other terms are not allowed as the metric has to be central. The inverse  is  
\[ (e^+,e^-)={1\over R_+(g_-)},\quad (e^-,e^+)={1\over R_-(g_+)}. \]
These functions correspond to square-lengths attached to the arrows as  $g_\pm(i)=g_{i\to i\pm 1}$ and it is natural to suppose that the metric is `edge symmetric' i.e. independent of the arrow direction\cite{Ma:sq}, which amounts to 
$g_-=R_-(g_+)$. In this case we have only one independent function $g:=g_+$ with $g_i$ being the square length of the link $i$--$i+1$,
\[ \cdots-\bullet_{i-1} {\buildrel g_{i-1}\over -}\bullet_i{\buildrel g_i\over -}\bullet_{i+1}{\buildrel g_{i+1}\over-}\bullet_{i+2}-\cdots.\]
So edge symmetric metrics correspond exactly to what you would expect for a metric on a lattice. We also let
\[ \rho_\pm=R_\pm({g_\pm\over g_\mp}),\quad \rho:=\rho_+\]
as  `ratio derivative' of the metric function $g$. 
Then \cite{ArgMa1} showed that there is a unique ($*$-preserving, left) QLC 
\[ \nabla e^\pm=(1-\rho_\pm) e^\pm\tens e^\pm,\quad \sigma(e^\pm\tens e^\pm)=\rho_\pm e^\pm\tens e^\pm ,\quad \sigma(e^\pm\tens e^\mp)=e^\mp\tens e^\pm.  \]
The exterior algebra here is supposed to be the standard one where $e^\pm$ are Grassmann and anticommute with each other. 

If we write $\nabla e^a=-\Gamma^a{}_{bc}e^b\tens e^c$ where indices run over $\pm$ then these same Christoffel symbols define $\nabla_\cX$ on a dual basis $f_a$ by 
\[ \nabla_\cX f_a=f_b\tens \Gamma^a{}_{bc} e^c.\]
 This applies when there is a basis of 1-forms, as is the case of any Cayley graph (where the vertices form a group and the arrows are right translation by a fixed set of generators). Now, in our case $\Gamma$ has only two nonzero entries namely
\[ \Gamma^\pm{}_{\pm\pm}=\rho_\pm-1\]
hence
\[ \nabla_\cX f_\pm = f_\pm\tens (\rho_\pm-1) e^\pm,\quad \sigma_\cX(e^\pm\tens f_\pm)= f_\pm\tens \rho_\pm e^\pm,\quad \sigma_\cX(e^\pm\tens f_\mp)= f_\mp\tens e^\pm.\]
Using the commutation rules,
\[ e_\pm f= R_\pm(f) e^\pm,\quad  f f_\pm= f_\pm R_\pm(f)\]
we check this is consistent with being a right connection,
\begin{align*} \nabla_\cX(f f_\pm)&=\sigma_\cX(\extd f\tens f_\pm)+ f f_\pm \tens (\rho_\pm-1) e^\pm=(\del_a f)\sigma_\cX(e^a\tens f_\pm)+ f_\pm \tens (\rho_\pm-1)e^\pm f\\
&=(\del_\mp f) f_\pm \tens e^\mp +  f_\pm \tens \rho_\pm e^\pm \del_\pm f+ f_\pm \tens (\rho_\pm-1) e^\pm f\\
&= f_\pm \tens (\del_\mp R_\pm f)e^\mp +f_\pm \tens e^\pm \del_\pm f+  f_\pm \tens (\rho_\pm-1) e^\pm \del_\pm f+ f_\pm \tens (\rho_\pm-1) e^\pm f\\
&=  f_\pm \tens (\del_a R_\pm f)e^a+ f_\pm \tens (\rho_\pm-1) e^\pm R_\pm f= \nabla_\cX(f_\pm R_\pm f)
\end{align*}
as required. Then
\[ \sigma_\cX^{-1}(f_\pm\tens e^\pm)={1\over R_\mp(\rho_\pm)}e^\pm\tens f_\pm,\quad  \sigma_\cX^{-1}(f_\pm\tens e^\mp)=e^\mp\tens f_\pm,\]
\begin{align*} {\rm div}(f_\pm)&=\ev(\sigma_\cX^{-1}\nabla_\cX f_\pm)=\sigma_\cX^{-1}(f_\pm(\rho_\pm-1)\tens e^\pm)=R_\mp(\rho_\pm-1)\ev(\sigma_\cX^{-1}(f_\pm\tens e^\pm))\\
&=R_\mp(\rho_\pm-1)\ev({1\over R_\mp(\rho_\pm)}e^\pm\tens f_\pm)= 1- { g_\mp\over g_\pm}\end{align*}
\begin{lemma} $\int f=\sum_i f\mu$ defined  by measure $\mu$ is divergence-compatible iff $\rho$ is  constant and ${R_\pm(\mu)\over\mu}=\rho^{\pm 1}$. Here, one may take $\mu=g$.
\end{lemma}
{\sl Proof.}  For the divergence compatibility, we need for all $f$, 
\[ \sum \mu f {\rm div} f_\pm = - \sum \mu f_\pm(\extd f)= -\sum\mu f_\pm((\del_a f)e^a)=-\sum\mu \del_\pm f=\sum\mu f-\sum R_\mp(\mu)f \]
which using $R_\pm(\mu)=\mu (1- {\rm div}(f_\mp))$, needs
\[ R_\pm(\mu)=\mu { g_\pm\over g_\mp};\quad R_+(\mu)=\mu {g\over R_- (g)},\quad R_-(\mu)=\mu {R_- (g)\over g}\]
which requires for a solution a constraint on $g$
\[ R_+(g) R_-(g)= g^2,\quad  g(i)=\left({ g(1)\over g(0)}\right)^i g(0). \]
We see that the condition on $g$ amounts to $\rho(i)=\rho(0)=g(1)/g(0)$, i.e. to $\rho$ a constant (this is equivalent geometrically to $\nabla$ having zero curvature) as well as 
\[ \rho_\pm=\rho^{\pm 1}\]
The simplest solution for $\mu$ is then $\mu=g$. $\hfill\square$\medskip

Proceeding for $\int$ divergence-compatible,  $e_\pm^*=-e_\mp$ implies
\[ f_\pm^*= -f_\mp {1\over R_{\mp}(\rho_\pm)}= - f_\mp {g_\mp\over g_\pm},\quad f_+^*=-{f_-\over\rho},\quad f_-^*=-\rho f_+.\]
Hence, writing $X=f_+ X^++ f_- X^-$, this is real with respect to our $*$ on vector fields iff
\begin{equation}\label{XrealZ} X^\pm{}^*=-\rho_\pm R_\pm(X^\mp)=-\rho^{\pm 1}R_\pm (X_\mp).\end{equation}
Next, from the above, we have 
\begin{equation}\label{kappaZ} \kappa={1\over 2}{\rm div}(X)={1\over 2}((\id- {R_-\over \rho})X^++ ((\id-\rho R_+)X^-),\end{equation}
which is real. We put this into the velocity equation (\ref{veleqX}), which for $\nabla_\cX$ above can be written as
\begin{align} \label{veleqZa} \dot X^+ &= \del_+(\kappa)X^+ + (1-\rho)X^+ X^+ -\rho\del_+(X^+)X^+- \del_-(X^+)X^- \\
 \label{veleqZb} \dot X^- &= \del_-(\kappa)X^- + (1-\rho^{-1})X^-X^--\rho^{-1}\del_-(X^-)X^-- \del_+(X^-)X^+\end{align}
where
\[ \del_+\kappa={1\over 2}((1-\rho R_+)\del_+ X^-+ (\del_++ {1\over\rho}\del_-)X^+),\quad 
\del_-\kappa={1\over 2}((1-{1\over \rho} R_-)\del_- X^++ (\del_-+\rho\del_+)X^-).\]
Now, applying $*$ to (\ref{veleqZb}) and subtracting from (\ref{veleqZa}) gives
\[ ((1-\rho)({1\over \rho^2}R_-- R_+) X^+)X^++((\rho\del_++\rho-1)\del_+ X^-)X^++(\del_- X^+)X^-+ \rho (\del_+ X^+)R_+^2(X^-)=0\]
which simplifies to the `improved auxiliary equation'
\begin{equation}\label{auxZ} (\del_-+\rho\del_+)(X^+R_+X^-)=(1-\rho)((R_+-{R_-\over \rho^2})X^+)X^+\end{equation}
where $X^+R_+X^-$ is real and $\del_-+\rho\del^+$ if $\rho=1$ would be the usual finite different Laplacian $\Delta_\Z$  (but is not the QRG Laplacian which would be  $\Delta =-{1+\rho\over h}\Delta_\Z$.) Note that if we take the * of the 2nd velocity equation and compare with the first, we obtain instead
\[ (\del_++{1\over \rho}\del_-)(X^-R_- X^+)=(1-{1\over\rho})(( R_--\rho^2 R_+)X^-)X^-\]
 which is equivalent to $*$ applied to (\ref{auxZ}). The latter when $\rho\ne 1$ amounts to the implicit condition
\[ (L X^+)X^+=\rho^2 R_+((L X^-)X^-);\quad  L=R_+-{1\over \rho^2}R_-.\]
and is needed precisely so that any initially real $X$ (in our sense (\ref{XrealZ})) stays real during the evolution, which in turn is needed for unitarity (\ref{unitarity}).

Once we have found the geodesic velocity field, it only remains to solve for the amplitude flow itself. This means
we solved for $\psi_s\in C(\Z)$ such that $\nabla_E\psi=0$, which is
\begin{equation}\label{ampZ} \dot \psi = -X_s(\extd \psi )-\psi \kappa_s= - (f_a X^a)((\del_b \psi)e^b)- \psi \kappa=-(\del_a\psi)X^a-\psi\kappa. \end{equation}

\subsection{Geodesics for the constant metric $\rho=1$ on $\Z$}

Here real $X$ means $X^+{}^*=-R_+(X^-)$ and $X^-{}^*=-R_-(X^+)$ and
\begin{align*} \kappa&=-{1\over 2}(\del_-X^++ \del_+X^-)\\
 2\del_+ \kappa&=-( R_+-1)(R_--1)X^+-\del_+^2 X^-=\Delta_\Z X^+ - \del_+^2 X^-.\end{align*}
The  geodesic velocity eqn is
\[ \dot X^+=(\del_+ \kappa)X^+-(\del_+ X^+)X^+-(\del_- X^+)X^-,\quad \dot X^-=(\del_-\kappa)X^--(\del_-X^-)X^--(\del_+X^-)X^+.\]
Applying $*$ to the second and comparing to the first gives 
\[ (\del_+^2X^-)X^++(\del_- X^+)X^-+(\del_+X^+)R_+^2 X^-=0\]
which simplifies to
\[ \Delta_\Z(X^+ R_+ X^-)=0,\]
where $X^+R_+X^-$ is real. Solving the first velocity equation and this `improved auxiliary equation' for real $X$ is equivalent to solving both velocity equations for real $X$. We first look at this latter moduli of solutions of the auxiliary equation as the subspace on which we will be solving the 1st velocity equation. Here among real non-negative functions
\[ \ker\Delta_\Z\cup\{{\rm non-negative}\}=\R_{\ge 0}1\]
(Here a general zero mode of $\Delta_\Z$ is $Y(i)=i \alpha - (i-1)\beta$ for $\alpha=Y(1),\beta=Y(0)$ but we need $\alpha\ge\beta\ge0$ for $Y$ to be positive for positive $i$ and $\beta\ge\alpha\ge0$ for $Y$ to be positive for  negative $i$.) We therefore can parameterize members of this aux moduli space as
\[  X^+=r e^{\imath \theta},\quad X^-=-R_-(re^{-\imath \theta});\quad r\in \R_{\ge 0},\quad \theta\in \C(\Z).\]

\begin{figure}
\[ \includegraphics[scale=.9]{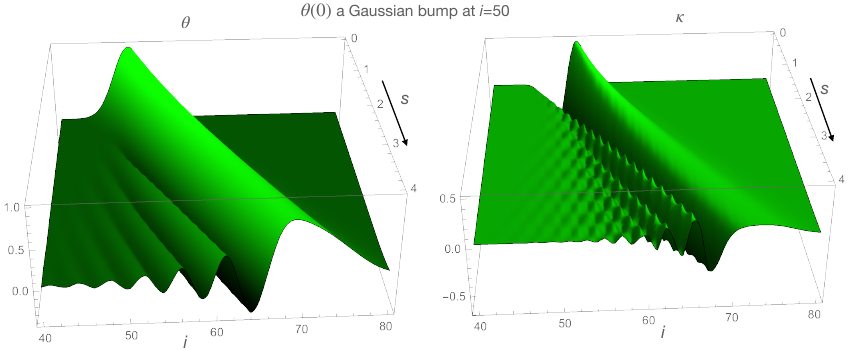}\]
\caption{(a) Smoothed view of solution of geodesic velocity equations for $\theta_s(i)$  with initial value at $s=0$ a Gaussian centred at $i=50$. Here $r=3$ resulting in an effective velocity of the peak of approximately 4.9 (b) Smoothed view of associated divergence $\kappa_s(i)$.\label{XsolZ}}
\end{figure}

All of this applies at each geodesic time $s$ and we now consider $r,\theta$  as functions of $s$ and solve the first vel equation which in terms of $X:=X^+,r$ looks like
\[ \dot X={1\over 2}(\Delta_\Z(X+ {r^2\over X}))X- (\del_+ X+ \del_-({r^2\over X}))X+ \del_- r^2={1\over 2}\left((R_--R_+)(X-{r^2\over X})\right)X\]
which, in terms of $r,\theta$, is
\[ \dot\theta= r(R_--R_+) (\sin(\theta)),
\quad \dot r=0\]
for the real and imaginary parts. Thus, $r$ becomes a constant parameter and any initial $\theta_0$ at $s=0$ evolves according to the equation shown for the angle phase $\theta$. For example, if we start at a sampled Gaussian peaked at $i=50$ then the flow is plotted (smoothly interpolated in $i$ for visual clarity) in Figure~\ref{XsolZ}(a). (By contrast, if we start at $\theta_0$ constant then this remains constant for all time.) The divergence for our form of $X^\pm$ is
\[ \kappa= -{1\over 2}\del_-(X+{r^2\over X})=-r\del_-\cos(\theta)\]
and plotted in Figure~\ref{XsolZ}(b).

\begin{figure}
\[ \includegraphics[scale=0.8]{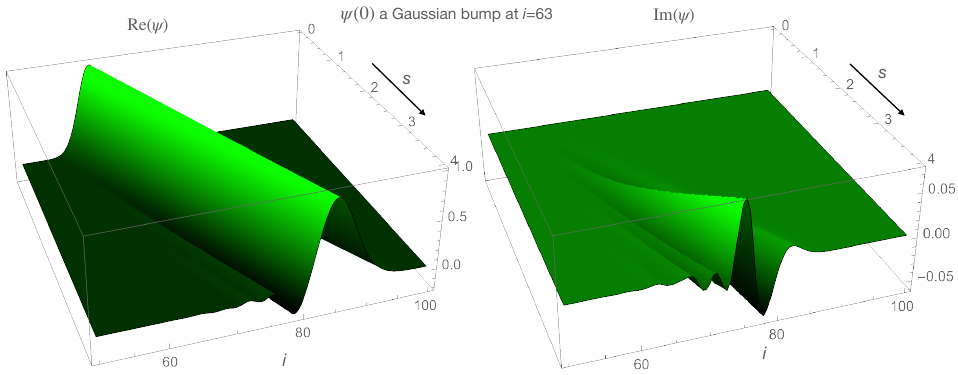}\]
\caption{Smoothed view of amplitude $\psi_s(i)$ for an initial Gaussian centred at $i=25$ and evolving under the geodesic velocity flow in Figure~\ref{XsolZ}. There is an effective velocity of approximately 5.6\label{psiZ} and the wave function acquires an imaginary component.}
\end{figure}

Finally, the amplitude flow is
\[ \dot\psi = - r(\del_+\psi)e^{\imath\theta}+ r(\del_-\psi)R_-( e^{-\imath\theta})-\psi \kappa.\]
An example is shown in Figure~\ref{psiZ}. We see that it acquires an imaginary component and becomes more wavelike with time. The imaginary component appears sooner if we start nearer (or on top of) the initial bump in $\theta$. Thus we see literally how a real bump function approximating a point in classical geometry evolves over geodesic time into a complex wave packet. 

\subsection{Geodesics for generic metrics on $\Z$}

For a generic metric on $\Z$, $\int$ defined by $\mu$ won't be divergence-compatible. But we can still define ${\rm div}_{\int}$ and require our vector field $X$ to be real with respect to $\int$.  Here
\[ \int X(\extd f)=\int (\del_\pm f) X(e^\pm)=\sum_i \mu (R_\pm f-f)X^\pm=\sum_i f (R_\mp (\mu X^\pm)-  \mu X^\pm)\]
(sum over the $\pm$ understood). For this to equal $-\int f {\rm div}_{\int}(X)$ for all $f$ as in (\ref{divint}), we need
\begin{equation}\label{divintZ} {\rm div}_{\int}(X)= -{1\over\mu}(\del_+(\mu X^+)+\del_-(\mu X^-)).\end{equation}
Next, $X$ is real with respect to $\int$ if the second of (\ref{unitarity}) holds, which means 
\[\int X((fe^\pm)^*)=\int X(e^\pm{}^* f^*)=-\int R_\mp(f^*) X(e^\mp )=-\int R_\mp(f^*) X^\mp=-\sum_i f^* R_\pm(\mu X^\mp)\]
has to equal
\[ \int X(fe^\pm)^*=\int f^*X^\pm{}^*=\sum_i f^* \mu X^\pm{}^*\] for all $f$. This gives us reality with respect to $\int$ as
\begin{equation}\label{XrealintZ} X^\pm{}^*=-{R_\pm(\mu X^\mp)\over\mu}.\end{equation}
We impose this on $X$ at all times $s$ and let $\kappa={1\over 2}{\rm div}_{\int}(X)$ so that the flow is unitary. One can check that $\kappa$ is then real. The velocity  equation (\ref{veleqX}) still gives  (\ref{veleqZa})-(\ref{veleqZb}) but with $\rho^{-1}$ in the latter replaced by $\rho_-$ since the metric is now arbitrary.  Comparison of one half of these velocity equations with the conjugate of the other half to give the improved auxiliary equation is now more complicated. One can also add a driving force term for maximal generality\cite{BegMa:gra}.

\end{document}